\documentclass{article}
\usepackage{amsfonts,amsmath,amstext,amsbsy,euscript,amssymb,amsthm}
 \font\amsy=msbm10
\def\IZ{\hbox{\amsy\char'132}}
\def\IR{\hbox{\amsy\char'122}}
\newtheorem{theorem}{Theorem}
\newtheorem{prop}{Proposition}
\newtheorem{definition}{Definition}
\begin{document}
\title{Black \& white Nim games}
\author{Urban Larsson}
\maketitle

\begin{abstract}
We present a new family of Nim games where the rules depend on a given `coloring' of the tokens, each token being either black or white. The rules are as in Nim with the restriction that a white token on top of each heap is not allowed. We resolve the winning strategies of two disjoint game families played on two heaps. The heap-sizes with black tokens correspond to the numbers $\lfloor \beta n \rfloor$, where $\beta>2$ is an integer for one of the families and irrational for the other, and where $n$ ranges over the positive integers. In the process, new notions of \emph{invariant games} are introduced.
\end{abstract}

\section{Introduction}
We assume that the reader is familiar with standard terminology on 2-player \emph{impartial games} \cite{BCG} and the game of Nim \cite{Bo}.
We present a new family of Nim games on a fixed finite number of heaps of `colored' tokens. Let $S\subset \IZ_{> 0}$. Then, for each heap, each $n^{th}$ token is \emph{black} if $n\in S$, otherwise \emph{white}. We say that a \emph{heap is black} if the top token of the heap is colored black \emph{or} if the heap is empty. Otherwise the \emph{heap is white}. The number of heaps stay constant for a given game, thus an empty heap is also ragarded as a heap.

For example, suppose that $S=\{2\}$ and there are three heaps with zero, two and three tokens respectively. That is the first heap is empty, the second heap has two tokens of which the lower is white and the upper black, the third has a black token in between two white ones. This means that the third heap is white, whereas the first and the second is black.  

Given any such coloring, the players move as in Nim on $k\in \IZ_{>0}$ piles, but where the \emph{legal positions} are restricted to that \emph{at least one} heap is black. That is a player may not move to a configuration where each heap is white. In our example this means that all Nim moves are allowed except the removal of precisely one token from the second heap.

The game constant $k$ determining the number of heaps follows the game until the end, irrespective of whether some heaps become empty (which they will). This means that if at least one heap is empty then all regular Nim-type moves are legal on the remaining heaps, \emph{independent of the coloring} and thus, the P-positions on the remaining heaps correspond precisely to those of Nim \cite{Bo}. Thus the interesting games are those where each heap contains a positive number of (colored) tokens.

For example, with $k=2$, if the heap sizes are 1 and 2 and only the second token (in the second heap) is colored black, then, the legal moves are $(1,2)\rightarrow (0,2)$ and $(1,2)\rightarrow (1,0)$. The second option is legal although `each' top token is white. Namely, the second heap is empty and therefore black.

We call the whole family of games \emph{black \& white Nim} and a particular game \emph{$k$-pile $S$-Nim}.

The games are \emph{symmetric}, by which we mean that the order of the heaps for a particular game is immaterial---and so the $k$-tuple $(x_1,\ldots ,x_k)$ is a P-position of $k$-pile $S$-Nim if and only if $(x_{\pi(1)},\ldots ,x_{\pi(k)})$ is a P-position for all permutations $\pi$ of $\{1,\ldots , k\}$. Henceforth we let the \emph{position} $(x_1,\ldots ,x_k)$ denote the unordered multiset
$\{x_1,\ldots ,x_k\}$, $x_i\in \IZ_{\ge 0}$, $i\in \{1,\ldots ,k\}$. %Let the set of legal positions, given a set $S$ and a $k$, be $\mathcal{S}\subset \IZ_{\ge 0}^k$ (that is a position is in $\mathcal{S}$ if and only if at least one heap is black). 

The \emph{legal moves} are independent of the position moved from, in the sense that, for all positions $X$ and $Y$, 
\begin{align}\label{XY}
X\rightarrow Y 
\end{align}
is a legal move if and only if it is a legal regular Nim move \emph{and} $Y$ is a legal position. Hence the moves in a given game can be regarded as \emph{invariant} in a similar sense as for the classical \emph{take-away} games \cite{Bo, W, G} extending the notion of invariant games from \cite{DR, LHF, L4, LW}. Therefore, let us give a new definition of game invariance for take-away games.

\subsection{Heap-size invariance}\label{S1.1} 
In some non-invariant (variant) take-away games, for example \cite{F,FL}, the players have to count the number of tokens in the piles in order to obtain the list of options. In our game this is not necessary. Many other take-away games \cite{L1,L2,L3,S,Z} are not variant in this sense, still they do not satisfy the criterion for [translation] invariance [of move options] defined in \cite{DR, LHF, L4}. (The brackets are ours.) The following notion of invariance in a broader sense is satisfied by the games in \cite{DR, LHF, L1,L2,L3,L4,S,Z} and many more. 

\begin{definition}
A take-away game on heaps of tokens is \emph{(heap-size) invariant} if the players do not need to count the number of tokens in the heaps in order to play by the rules.
\end{definition}

\begin{prop}\label{P}
Our black \& white Nim games are heap-size invariant.
\end{prop}

\noindent{\bf Proof.} 
Let $S\subset \IZ_{>0}$ and $k\in \IZ_{>0}$. Suppose that $X=(x_1,\ldots ,x_k)$ is a (legal) non-terminal position of $k$-pile $S$-Nim. Then each Nim move $x_i \rightarrow y$, $i\in \{1,\ldots , k\}$ is legal or illegal irrespective of the size of the heaps. The move is legal if and only if at least one of the heaps in the new position $(x_1,\ldots,x_{i-1},y ,x_{i+1},\ldots, x_k)$ is black. The coloring of the tokens is done before the game starts and hence, to play by the rules, specific information about the number of tokens in the heaps is not needed. \hfill $\Box$\\

In the coming, we assume that the reader is acquainted with the notion of invariance used in \cite{DR,LHF}, otherwise jump to Section \ref{S2}.

\subsection{$(\mathcal{M}, \mathcal{B})$-invariance}
In view of Section \ref{S1.1}, an alternative definition of invariance which extends that in \cite{DR,LHF,L4} follows. We assume that the reader is familiar with the notation used in \cite{LHF} and, of course, as before, the meaning of invariance in that paper.

\begin{definition}\label{D}
Let $k\in \IZ_{>0}$, $S\subset \IZ_{>0}$, $\mathcal{M}\subset (\IZ_{\ge 0} ^{\;k}\setminus \{\boldsymbol 0\})$ and $\mathcal{B}\subset \IZ_{\ge 0}$. The legal positions belong to the set $\mathcal{B}$. The moves belong to the set $\mathcal{M}$. In actual play, a move $\boldsymbol m\in \mathcal{M}$ from $\boldsymbol x\in \mathcal{B}$ is legal if and only if $(\boldsymbol x\ominus \boldsymbol m)\in \mathcal{B}$. A game which satisfies this form of invariance is called $(\mathcal{M},\mathcal{B})$-invariant.
\end{definition}

\begin{prop}
Let $S\subset \IZ_{>0}$ and $k\in \IZ_{>0}$. Our black \& white Nim games are $(\mathcal{M},\mathcal{B})$-invariant where $\mathcal{M}$ denotes the set of all Nim moves (that is $k$-tuples of the form $(x,0,\ldots,0)$, $x\in \IZ_{>0}$, with unordered notation as before) and where  $\mathcal{B} = (S\cup\{0\})^k\subset \IZ_{\ge 0}$.
\end{prop}
 
We omit the proof which is by construction. 

Observe that $(\mathcal{M}, \IZ_{\ge 0} ^{\;k})$-invariance corresponds to the notion of invariance used in \cite{DR,LHF,L4}. Any $(\mathcal{M},\mathcal{B})$-invariant game is of course heap-size invariant. This follows by a similar argument as in the proof of Proposition \ref{P}.

\section{Modular and Beatty black \& white Nim}\label{S2}
Notice that, by the rules of regular Nim, the $X$ in (\ref{XY}) above can be any position, that is, we do not have to require for a \emph{starting position} of black \& white Nim to be legal. If a non-terminal position is not a starting position then, by the rules of Nim, there is a legal move, for example, remove all tokens from one of the non-empty heaps. For simplicity, in the remainder of this note, we require that also the starting position contains at least one black heap.

Let $\beta\in \IR$, $\beta \ge 2$. We resolve black \& white Nim on two heaps for the set $S = \{\lfloor \beta n \rfloor \mid n\in \IZ_{>0} \}$ and where $\beta$ is either irrational or an integer. 

The second result, which concerns the irrational case, inspired this note in the first place. We call those games \emph{Beatty Nim} since a sequence of the form $(\lfloor n\beta \rfloor)$, for a positive irrational $\beta$ and where $n$ ranges over the positive integers, is known as a \emph{Beatty sequence}. Two sequences of positive integers are said to be \emph{complementary} if each positive integer occurs precisely once in precisely one of these sequences. The celebrated Beatty's theorem \cite{Be} (discovered by lord Rayleigh \cite {R}) states that a pair of Beatty sequences $(\lfloor \alpha n \rfloor)$ and $(\lfloor \beta n \rfloor)$ are complementary if and only if $\alpha^{-1}+\beta^{-1}=1$ and $\beta$ is a positive irrational.

A conjecture in \cite{DR}, which is resolved in \cite{LHF}, asks for games with invariant move sets, for which complementary pairs of Beatty sequences constitute the set of P-positions. The resolution in \cite{LHF} is satisfactory in a sense, namely the game rules produced by that result are invariant (in both senses). On the other hand, the game rules seem very complex, which is unusual for combinatorial games, where game rules are often very simple, allowing for a child to easily understand the rules. 

At the GONC 2011 conference A. S. Fraenkel asked for `short game rules' given a specific pair of complementary Beatty sequences. In this paper we study a game family which satisfies both requirements, short rules and invariance, the hinge of course being the pre-coloring of the tokens. But, as soon as the tokens have been colored (which can be done in $\log(n)$ polynomial time by a computer) the game rules are nearly as easy as those of regular Nim. 

The first result concerns a related family of games which we denote by \emph{modular black \& white Nim}.

\begin{theorem}\label{one}
Let $\beta\ge 2$ denote an integer and let $S = \{\beta n \mid n\in \IZ_{>0}\}$. Then a position $(x,y)$ of 2-pile $S$-Nim is a previous player winning position if and only if 
\begin{align}\label{xy}
(x,y) = (\beta (n+t), \beta n + t), 
\end{align}
for some $n\in \IZ_{>0}$ and where $t\in \{1,\ldots , \beta -1\}$. 
\end{theorem}

\noindent {\bf Proof.} Suppose that  $x = \beta (n+t)$ and $y=\beta n + t$, for some $n\in \IZ_{>0}$ and some $t\in \{1,\ldots , \beta -1\}$. Then, we have to show that none of the options from $(x,y)$ is of the same form. If a player moves on the y-heap, then if the congruence class changes, the new position is clearly not of the same form, so suppose that $y\rightarrow \beta m +t$, for some integer $0\le m<n$. Then, the new position is $(\beta (n+t), \beta m +t)$ with $m<n$, which is not of the form in (\ref{xy}).

Suppose on the other hand that $(x,y)$ is not of the form in (\ref{xy}). Then we have to show that there is an option to a position of this form. By the definition of the set $S$, at least one position is of the form $\beta i$, $i\in \IZ_{\ge 0}$. Hence suppose that $x = \beta i$, $i\in \IZ_{\ge 0}$. If $y\ge x$, then there is a lowering of the y-heap to a position of the form in (\ref{xy}), for example $(\beta i, \beta (i-1)+1)$ will do. 

Hence suppose that $\beta i = x > y = \beta n +t$, for some $n\in \IZ_{\ge 0}$ with $t\in \{1,\ldots ,\beta-1\}$. (If $y$ were of the same form as $x$ we would swap piles and return to the previous paragraph hence we exclude the case $t=0$ here.) This gives that $i> n+t/\beta$ which implies that $i \ge n+1$. We have already dealt with the case $i= n+t$ in the first part of the proof. Suppose first that $i>n+t$. Then a lowering of the x-heap to $n+t$ is the desired move. If $n<n+(i-n)<n+t$, then a lowering of the y-heap to $\beta n + (i-n)$ suffices. \hfill $\Box$\\

%\noindent {\sc Case 1:} Suppose that $i>n+t$. Then a lowering of the x-heap to $\beta (n+t)$ suffices.\\

%\noindent {\sc Case 2:} Suppose that $i<n+t$. Then a lowering of the y-heap to $\beta (i-t)+t$ suffices.

\begin{theorem}\label{two}
Let $\beta>2$ denote an irrational and let $S = \{\lfloor \beta n \rfloor \mid n\in \IZ_{>0}\}$. Then a position $(x,y)$ of 2-pile $S$-Nim is a previous player winning position if and only if $(x,y)=(\lfloor \alpha n \rfloor, \lfloor\beta n\rfloor)$ for some $n\in \IZ_{\ge 0}$ and where 
\begin{align}\label{Beatty}
\frac{1}{\alpha} + \frac{1}{\beta} = 1.
\end{align}
\end{theorem}

\noindent {\bf Proof.} Let the position be $(x,y)$. 

Suppose at first that $x=\lfloor \alpha n \rfloor$ and $y = \lfloor\beta n\rfloor$, for some $n\in \IZ_{\ge 0}$. If $n=0$, we are done, so suppose that $n>0$. Observe that any position of this form is legal since one of the coordinates is $\lfloor \beta n\rfloor$. Then we have to show that there is no Nim option of the same form. But this follows by (\ref{Beatty}) which, by Beatty's theorem, implies that the sequences $(\alpha i)$ and $(\beta i)$ are complementary. 

Otherwise, there are non-negative integers $m\ge n$ such that either\\

\noindent{\sc Case 1:} $x=\lfloor \alpha m \rfloor, y= \lfloor\beta n\rfloor$, $m>n$\\

\noindent{\sc Case 2:} $x=\lfloor \alpha m \rfloor, y=\lfloor\alpha n\rfloor$\\

\noindent{\sc Case 3:} $x=\lfloor \alpha n \rfloor, y=\lfloor\beta m\rfloor$, $m>n$\\

\noindent{\sc Case 4:} $x=\lfloor \beta m \rfloor, y=\lfloor\beta n\rfloor$, $n>0$.\\

Let us begin by excluding the position $(x,y)$ given by the second case. It is not legal. This follows from complementarity of the sequences $(\lfloor\alpha i\rfloor)$ and $(\lfloor\beta i\rfloor)$, namely since $x = \lfloor \alpha n \rfloor$ there is no integer $i$ such that $\lfloor \beta i\rfloor = x$ and similarly for $y$.

The three remaining cases represent legal positions. Notice that none represents a position of the form in the theorem, the first and third since the sequences are strictly increasing and the fourth by complementarity. Hence, our task is to find a legal move to a position of the form in the theorem, for each case. 
For the first case, since $m>n$ we can lower the x-heap to $\lfloor\alpha n\rfloor$, which gives a position of the desired form. For the third case, by $m>n$ since $\beta>2$ we get $\lfloor \beta m\rfloor>\lfloor \beta n\rfloor$, so that the desired Nim move on the y-heap is to lower it to the position $(\lfloor \alpha n \rfloor, \lfloor\beta n\rfloor)$. The fourth case is similar, but the lowering is on the x-heap, motivated by $\lfloor \beta m \rfloor\ge \lfloor\beta n\rfloor >\lfloor \alpha n\rfloor $, which follows since (\ref{Beatty}) gives $1<\alpha<2<\beta $ and by $n>0$. (The latter inequality excludes the terminal position $(x,y) = (0,0)$ which of course is also of the form $(x,y) = (\lfloor \alpha n \rfloor, \lfloor\beta n\rfloor)$, for some $n\in \IZ_{\ge 0}$).
\hfill $\Box$

\section{Remarks---open problems}
In view of Theorem \ref{one} and \ref{two}, what are the P-positions of the 2-pile take away games whenever $\beta>2$ is a non-integer rational?

Given a set $S$, one can generalize our games by allowing any move set, as in \cite{DR, G, LHF, L4} (not only Nim type moves), but keep the rest of the heap rules as defined in the first section. 

Another generalization is to allow the piles to have different black \& white colorings (and maybe start looking at only Nim type moves). Even more generally, allow more colors and put various chromatic restrictions on the type of legal moves. 

For example, in the games of $l$-colored $k$-pile Nim, let each token take one of $l$ colors. The players move as in Nim. One game in this family is the game of \emph{Spectrum Nim}: Let the tokens have colors of a modulo $l$ \emph{spectrum}, so that the $n^{th}$ token (in each pile) has color $i$, where $i\equiv n\pmod l$. We suggest two variations of the game:\\

\noindent{\sc (Full) Spectrum Nim:} At each stage of the game, all colors must be represented in case $l\ge k$. Otherwise all non-empty heaps have different colors.\\

\noindent{\sc Bi-chromatic Nim:} This game has three variations, at each stage of the game, at most/precisely/at least two distinct colors must be represented among the $k$ heaps.\\

What are the P-positions (Grundy values) for the respective games? 

For a variation of black \& white Nim, one could study partizan variations of the games in Theorem 1 and 2 where Left has to move to a position with at least one white heap and Right has to move to a position with at least one black heap.\\

\noindent{\bf Acknowledgments.} I would like to thank Aviezri S. Fraenkel and Mike Weimerskirch for many interesting discussions.


\begin{thebibliography}{11}
\bibitem[BCG]{BCG} E. R. Berlekamp, J. H. Conway, R. K. Guy,  \emph{Winning ways}, {\bf 1-2} Academic Press, London (1982). Second edition, {\bf 1-4}. A. K. Peters, Wellesley/MA (2001/03/03/04).
\bibitem[Be]{Be} S. Beatty, Problem 3173, \emph{Amer. Math. Monthly}, {\bf 33} (1926) 159.%, {\bf 34} (1927) 159-160.
\bibitem[Bo]{Bo} C. L. Bouton, Nim, A Game with a Complete Mathematical Theory
\emph{The Annals of Mathematics}, 2nd Ser., Vol. 3, No. 1/4. (1901 - 1902), pp. 35-39.
\bibitem[DR]{DR} E. Duch\^{e}ne and M. Rigo, Invariant Games,
\emph{Theoret. Comp. Sci.}, Vol. 411, 34-36 (2010), pp. 3169-3180 
%\bibitem[O'B03]{OB2003} K. O'Bryant, Fraenkel's Partition and Brown's
%Decomposition \emph{Integers}, {\bf 3} (2003), A11, 17 pp.
%\bibitem[Fr69]{F1969} A.S. Fraenkel, The bracket function and complementary
%sets of integers, \emph{Canad. J. Math.} {\bf 21} (1969), 6-27.
%\bibitem[Fr82]{F1982}
%A.\ S.\ Fraenkel, How to beat your Wythoff games' opponent on
%three fronts, {\it Amer. Math. Monthly\/} {\bf 89} (1982) 353--361.
%\bibitem[Fr08]{F2010} A.\ S.\ Fraenkel, The Rat game and the Mouse game, 
%to appear in \emph{Games of no Chance 2008}.
%\bibitem [F2]{F2} A.\ S.\ Fraenkel, Complexity, appeal and challenges 
%of combinatorial games, Proc. of Dagstuhl Seminar ``Algorithmic 
%Combinatorial Game Theory'', \emph{Theoret. Comp. Sci.} 
%{\bf 313} (2004) 393-415, special issue on Algorithmic Combinatorial Game Theory.
\bibitem[F]{F} A.\ S.\ Fraenkel, The Rat game and the Mouse game, to appear in \emph{Games of no Chance 2008}.
\bibitem[FL]{FL} A.\ S.\ Fraenkel, U.\ Larsson, A variant Wythoff Nim, preprint.
\bibitem[G]{G} S. W.\ Golomb, A mathematical investigation of games of "take-away''. \emph{J. Combinatorial Theory} {\bf 1} (1966) 443---458.
\bibitem[L1]{L1} U.\ Larsson Restrictions of $m$-Wythoff Nim and $p$-complementary Beatty sequences, to appear in \emph{Games of no Chance 2008}.
\bibitem[L2]{L2} U.\ Larsson, 2-pile Nim with a Restricted Number of Move-size Imitations, \emph{Integers} {\bf 9} (2009), Paper G4, pp 671-690.
\bibitem[L3]{L3} U.\ Larsson, Blocking Wythoff Nim, \emph{The Electronic Journal of Combinatorics}, P120 of Volume 18(1) (2011).
\bibitem[L4]{L4} U.\ Larsson, The star operator and invariant subtraction games, preprint.
\bibitem[LHF]{LHF} U.\ Larsson, P.\ Hegarty, A.\ S.\ Fraenkel, Invariant and dual subtraction games resolving the Duch\^ene-Rigo Conjecture, \emph{Theoret. Comp. Sci.} Vol. 412, 8-10 (2011) pp. 729-735.
\bibitem[LW]{LW} U.\ Larsson, M.\ Weimerskirch, Impartial games, whose rule sets correspond to a given continued fraction, preprint.  
\bibitem[R]{R}J. W. Rayleigh. The Theory of Sound, \emph{Macmillan, London}, (1894) p. 122-123.
%\bibitem[L3]{L3} Permutation games and the $\star$-operator, U. Larsson, 
%preprint.
%\bibitem[LW]{LW} U. Larsson, J. W\"{a}stlund, From heaps of matches to 
%undecidability of games, preprint.
%\bibitem[HL06]{HL2006} P. Hegarty and U. Larsson, Permutations of the 
%natural numbers with prescribed 
%difference multisets, \emph{Integers} {\bf 6} (2006), Paper A3, 25pp. 
\bibitem[S]{S} A. J. Schwenk, ``Take-Away Games", \emph{Fibonacci Quart.} {\bf 8} (1970), 225-234.
\bibitem[Z]{Z} Michael Zieve, Take-Away Games, \emph{Games of No Chance, MSRI Publications}, {\bf 29}, (1996) pp. 351–361
\bibitem[W]{W} W.A. Wythoff, A modification of the game of Nim, \emph{Nieuw Arch. Wisk.} {\bf 7} (1907) 199-202.
\end{thebibliography}
\end{document}